\documentclass{article}
\usepackage{graphicx,newpxmath,newpxtext,mathrsfs} 

\title{Martingale expansion for stochastic volatility}
\author{Masaaki Fukasawa\\
{\small Graduate School of Engineering Science}\\
{\small The University of Osaka}\\
{\small 560-8531 Japan}}
\date{}

\newtheorem{thm}{Theorem}
\newtheorem{cor}{Corollary}
\newtheorem{rem}{Remark}
\begin{document}
\maketitle
\begin{abstract}
  The martingale expansion provides a refined approximation to the marginal distributions of martingales beyond the normal approximation implied by the martingale central limit theorem. We develop a martingale expansion framework specifically suited to continuous stochastic volatility models. Our approach accommodates both small volatility‑of‑volatility and fast mean‑reversion models, yielding first‑order perturbation expansions under essentially minimal conditions.
\end{abstract}

\section{Introduction}
Stochastic volatility (SV) models constitute a central class of continuous‑time asset price models in which the instantaneous variance of returns is itself governed by an additional stochastic process. Unlike the classical Black–Scholes framework, which assumes constant volatility, SV models allow volatility to evolve randomly over time, thereby capturing a range of empirically observed features of financial markets, including volatility clustering, heavy‑tailed return distributions, and pronounced implied‑volatility smiles and skews.

Formally, an SV model specifies the joint dynamics of an asset price process
$S$ and its latent spot variance process $V$, typically through a system of coupled stochastic differential equations. The spot variance process $V$ is often mean‑reverting and may be correlated with the asset‑price shocks, a feature that enables the model to reproduce leverage effects observed in equity markets. Prominent examples include the Heston model, in which $V$ follows a square‑root diffusion, and the SABR model, widely used in interest‑rate markets due to a closed form approximation formula for the implied volatility.

The flexibility afforded by stochastic volatility has made SV models indispensable in modern derivative pricing and risk management. They provide a more realistic representation of market dynamics and yield option prices that align more closely with observed implied‑volatility surfaces. At the same time, the introduction of a latent volatility factor complicates both analytical tractability and statistical inference, motivating a substantial literature on approximation methods, asymptotic expansions, and efficient numerical schemes.

Consider an abstract SV model with zero interest rates
\begin{equation}\label{model}
 \frac{\mathrm{d}S_t}{S_t} = \sqrt{V_t}\mathrm{d}B_t, \ \ 
 B_t=\rho W_t+ \sqrt{1-\rho^2}W^\perp_t
\end{equation}
for an asset price process $S$,
where $(W,W^\perp)$ is a $2$-dimensional standard Brownian motion  on a filtered probability space $(\Omega,\mathscr{F},\mathsf{P}, \{\mathscr{F}_t\})$,
$\rho\in (-1,1)$,
and $V$ is a nonnegative cadlag process adapted to a smaller filtration $\{\mathscr{G}_t\}$
to which $W$ is also adapted while $W^\perp$ is independent.
Suppose that $V = V^\epsilon$ depends on a parameter $\epsilon > 0$ and
as $\epsilon \to 0$, 
\begin{equation*}
    \int_0^T V^\epsilon_t \mathrm{d}t - v^\epsilon
\end{equation*}
converges to $0$ in probability for a deterministic positive sequence $v^\epsilon$ with limit
\begin{equation}\label{liminf}
  v^0:=  \lim_{\epsilon \to 0} v^\epsilon >0.
\end{equation}
Then $S = S^\epsilon$ also depends on $\epsilon$ and 
by the martingale central limit theorem, 
$\log S_T = \log S^\epsilon_T$ converges in law to the normal distribution with mean 
$-v^0/2$ and variance $v^0$.
The model \eqref{model} with small $\epsilon >0$ is  close to the Black-Scholes model in this sense.

Since SV models generally do not admit closed-form expressions for derivative prices or hedging strategies, asymptotic expansions 
with respect to $\epsilon$ for various models of $V^\epsilon$
have been extensively developed for both practical implementation and theoretical analysis. Examples include
small volatility-of-volatility models, e.g.~\cite{Lewis,BG, Alos, Gobet, BFG}, where $\epsilon$ is the volatility parameter of $V^\epsilon$, and fast-mean-reverting models, e.g.~\cite{Fouque,FukasawaE,GS},
where $V^\epsilon$ is ergodic and a certain negative power of $\epsilon$ is the time scale parameter of $V^\epsilon$.

In \cite{FukasawaM}, the author introduced a unified framework for computing and validating first-order asymptotic expansions, based on Yoshida’s martingale expansion theory \cite{Yoshida1,Yoshida2} combined with a partial Malliavin calculus.
In the present paper, we propose an alternative approach that is more direct and elementary, and which in particular establishes the validity of the first-order expansion under weaker and essentially minimal conditions.

The martingale expansion was first formulated by Mykland~\cite{Myk1,Myk2} for twice continuously differentiable test functions using It\^o's formula, and then by Yoshida~\cite{Yoshida1,Yoshida2} for a general test function under a non-degeneracy condition on Malliavin covariances.
In \cite{FukasawaM}, we relied on the condition
$|\rho|<1$ that provides a smoothness of the distribution of $S_T$, admitting an effective application of the partial Malliavin calculus to ensure the required non-degeneracy.
In the present paper, we directly utilize a smoothness property due to $|\rho|<1$  to bypass the Malliavin calculus.

\section{Results}
Here we state the main results of the paper. All the proofs are deferred to Section~\ref{proof}.
Let $T>0$ be fixed and
\begin{equation*}
    (X^\epsilon, Y^\epsilon) =
    \left(\frac{1}{\sqrt{v^\epsilon}}\int_0^T \frac{\mathrm{d}S^\epsilon_t}{S^\epsilon_t},
    \frac{1}{\epsilon}\left(\int_0^T V^\epsilon_t \mathrm{d}t - v^\epsilon\right) \right).
\end{equation*}
Let $\phi$ denotes the standard normal density function.
\begin{thm}\label{thm}
If $Y^\epsilon$ is uniformly integrable and if
$(X^\epsilon,Y^\epsilon)$ converges in law 
as $\epsilon \to 0$ to such $(X,Y)$ 
that $x \mapsto \mathsf{E}[Y|X=x]$ is twice continuously differentiable 
with
\begin{equation}\label{boundary}
    \lim_{|x|\to \infty}\mathsf{E}[Y|X = x] \phi(x) =  \lim_{|x|\to \infty} \frac{\mathrm{d}}{\mathrm{d}x} (\mathsf{E}[Y|X = x] \phi(x)) = 0
\end{equation}
and that
\begin{equation*}
        \phi^\epsilon(x) := \phi(x) + \frac{\epsilon}{2\sqrt{v^\epsilon}}\frac{\mathrm{d}}{\mathrm{d}x}
        (\mathsf{E}[Y|X = x] \phi(x))
         + \frac{\epsilon}{2v^\epsilon}\frac{\mathrm{d}^2}{\mathrm{d}x^2}
        (\mathsf{E}[Y|X = x] \phi(x))
    \end{equation*}
is integrable on $\mathbb{R}$,
    then for any bounded Borel function $f$,
    \begin{equation}\label{main}
        \mathsf{E}[f(S^\epsilon_T)] = \int f\left(S_0 \exp\left\{\sqrt{v^\epsilon}x - \frac{v^\epsilon}{2} \right\}\right)\phi^\epsilon(x)\mathrm{d}x + o(\epsilon).
    \end{equation}
\end{thm}
\begin{rem}\label{rem1}
    \upshape
    By the martingale central limit theorem, the martingale marginal $X^\epsilon$ converges to the standard normal distribution if $\epsilon Y^\epsilon \to 0$. Therefore,
    $X$ is always a standard normal random variable under the assumption.
As illustrated in \cite{FukasawaM},  for both small volatility-of-volatility and fast-mean-reverting models, by suitably choosing $v^\epsilon$, it is straight-forward to observe that $(X^\epsilon,Y^\epsilon)$ converges in law to a  $2$-dimensional centered normal distribution. In such a case, $Y$ is equal in law to $\mathsf{E}[XY]X + Z$ with a centered normal random variable $Z$ with $\mathsf{E}[XZ] = 0$. In particular, $\mathsf{E}[Y|X=x] = \mathsf{E}[XY]x$, and $\eqref{boundary}$ is then trivial.
\end{rem}
\begin{cor}\label{cor1}
Under the assumption of Theorem~\ref{thm}, for any $K>0$,
\begin{equation}\label{put}
\begin{split}
\mathsf{E}[(K-S^\epsilon_T)_+] &= p_K(S_0,v^\epsilon ) + \frac{\epsilon}{2\sqrt{v^\epsilon}} K \mathsf{E}[Y|X=-d_-(S_0,v^\epsilon )] \phi(-d_-(S_0,v^\epsilon )) +o(\epsilon)\\
&=p_K(S_0,v^\epsilon ) + \epsilon \frac{\partial p_K}{\partial t}(S_0,v^\epsilon)\mathsf{E}[Y|X=-d_-(S_0,v^\epsilon )] +o(\epsilon)\\
&=p_K\left(S_0,v^\epsilon + \epsilon\mathsf{E}[Y|X=-d_-(S_0,v^\epsilon )] +o(\epsilon)\right) ,
\end{split}
\end{equation}
where
\begin{equation}\label{bsput}
    \begin{split}
        &p_K(s,t) = K\Phi(-d_-(s,t)) -s \Phi(-d_+(s,t)),\\
        &d_{\pm}(s,t) = \frac{1}{\sqrt{t}} \left(\log \frac{s}{K} \pm \frac{t}{2} \right)
    \end{split}
\end{equation}
and $\Phi$ is the standard normal distribution function.
\end{cor}
\begin{rem}
    \upshape
    Note that $p_K(s,t)$ in \eqref{bsput} is the Black-Scholes put option price function with total variance $t$.
    The last expression of \eqref{put} implies that the Black-Scholes implied total variance $\hat{v}(K)$, which is defined by
     \begin{equation}\label{imp}
        \mathsf{E}[(K-S^\epsilon_T)_+] = p_K(S_0,\hat{v}(K) ),
    \end{equation}
    is expanded as
    \begin{equation}\label{ivex}
        \hat{v}(K) =v^\epsilon + \epsilon\mathsf{E}[Y|X=-d_-(S_0,v^\epsilon )] + o(\epsilon).
    \end{equation}
    Although this is valid for any $K>0$ as an asymptotic formula, taking into account how it is derived, we can expect its numerical accuracy only when
    \begin{equation*}
        \frac{\partial p_K}{\partial t}(S_0,v^\epsilon) = \frac{K}{2\sqrt{v^\epsilon}} \phi(-d_-(S_0,v^\epsilon)) 
        = \frac{S_0}{2\sqrt{v^\epsilon}} \phi(-d_+(S_0,v^\epsilon)) 
    \end{equation*}
    is not too small, which is the case $K$ is near $S_0$, i.e., around the at-the-money.
    Since $(X,Y)$ is the limit of $(X^\epsilon,Y^\epsilon$), and so, is the limit of
    \begin{equation*}
        \left(\frac{1}{\sqrt{v^\epsilon}} \left( \log \frac{S_T}{S_0}  + \frac{v^\epsilon}{2} \right), Y^\epsilon \right),
    \end{equation*}
    by formally replacing $(X,Y)$ with the above, we reach a conceptually interesting formula
    \begin{equation*}
        \hat{v}(K) \approx \mathsf{E}\left[ \int_0^T V^\epsilon_u \mathrm{d}u \bigg| S_T = K \right].
    \end{equation*}
    A rigorous validation of this formula is left for future research.
\end{rem}

\begin{rem}
    \upshape
    As is well-known,
    differentiating the defining equation \eqref{imp} of the implied total variance in log moneyness,
    we reach a volatility skew formula
    \begin{equation*}
        \frac{\partial}{\partial k} \sqrt{\hat{v}(S_0e^k)}  = \frac{\mathsf{P}[S_0e^k > S^\epsilon_T] - \Phi(-d_-(S_0,\hat{v}(S_0e^k)))}{\phi(-d_-(S_0,\hat{v}(S_0e^k)))}.
    \end{equation*}
Here, the derivative exists because under the current assumption of $|\rho| < 1$, $\log S^\epsilon_T$ follows a mixed normal distribution and in particular admits a density.
We have asymptotic expansions for both $\mathsf{P}[K > S^\epsilon_T] $ and 
$\hat{v}(K)$ from \eqref{main} with $f(s) = 1_{(-\infty,K)}(s)$ and \eqref{ivex}. Substituting those, we obtain
 \begin{equation*}
        \frac{\partial}{\partial k} \sqrt{\hat{v}(S_0e^k)}  = \frac{\epsilon}{2\sqrt{v^\epsilon}} \frac{\partial}{\partial k}
        \mathsf{E}\left[
        Y\bigg| X = \frac{k}{\sqrt{v^\epsilon}} + \frac{\sqrt{v^\epsilon}}{2}
        \right] + o(\epsilon).
    \end{equation*}
    In the case $\mathsf{E}[Y|X=x] = \mathsf{E}[XY]x$ mentioned in Remark~\ref{rem1}, we have
    \begin{equation}\label{skew}
        \frac{\partial}{\partial k} \sqrt{\hat{v}(S_0e^k)}  = \frac{\epsilon}{2v^\epsilon}
        \mathsf{E}\left[
        XY
        \right] + o(\epsilon).
    \end{equation}
\end{rem}
\begin{rem}
    \upshape
    In \cite{FukasawaM}, a more abstract framework is given, which in particular incorporates jumps and time-dependent $\rho$. For a continuous model~\eqref{model}, it however requires a stronger  moment conditions on $(X^\epsilon,Y^\epsilon)$.
\end{rem}
\section{Example}
In \cite{FukasawaM}, we have already observed how the martingale expansion framework accommodates various perturbation models including small volatility-of-volatility and fast-mean-reverting models. Here, for an illustrative purpose, we take the small volatility-of-volatility expansion of a Bergomi-type model as an example.

Consider
\begin{equation*}
    V^\epsilon_t = V_0(t)\exp\left\{
\epsilon \sum_{i=1}^d \int_0^t k_i(t-s)\mathrm{d}W^i_s - \frac{\epsilon^2}{2}
\sum_{i=1}^d \int_0^t k_i(t-s)^2 \mathrm{d}s
    \right\}
\end{equation*}
with a small volatility-of-volatility parameter $\epsilon>0$,
where $V_0(t)$ is a deterministic positive continuous function,
$(W^1,\dots,W^d)$ is a $d$-dimensional standard Brownian motion correlated with $B$ in \eqref{model} as 
\begin{equation*}
    \mathrm{d}\langle B, W^i\rangle_t = \rho_i \mathrm{d}t,  \ \  
    \rho := \sqrt{\sum_{i=1}^d \rho_i^2}  < -1,
\end{equation*}
 and $k_i$ are locally square integrable functions on $[0,\infty)$.

The function $t \mapsto V_0(t)$ is called the forward variance curve (at time $0$), due to
$V_0(t) = \mathsf{E}[V^\epsilon_t]$.
The case $k_i(t) = a_ie^{-b_it}$, $a_i, b_i >0$, 
describes the Bergomi model (see \cite{Bergomi}).
The case $d=1$ with $k_1(t) = a t^{H-1/2}$, $a>0$, $H \in (0,1/2)$, describes the rough Bergomi model proposed by \cite{BFG}.
In \cite{EFGR,FukasawaV, FukasawaW}, 
short-time (non-small volatility-of-volatility) expansions of the implied volatility and skew for this type of model are given.
We refer the reader to Section~8.3 of \cite{BFFGJR} for the difference between short-time and small volatility-of-volatility expansions.

Consistently to \eqref{model}, we have a decomposition
\begin{equation*}
    B =  \rho W + \sqrt{1-\rho^2} W^\perp, \ \ 
    W = \frac{1}{\rho} \sum_{i=1}^d \rho_iW^i, \ \ 
    \ \ W^\perp = \frac{1}{\sqrt{1-\rho^2}}\left(B -
     \sum_{i=1}^d \rho_iW^i\right)
\end{equation*}
with $W^\perp$ being a standard Brownian motion independent of $(W^1,\dots,W^d)$.
We take as $\{\mathscr{G}_t\}$ the natural filtration generated by $(W^1,\dots,W^d)$.

By taking
\begin{equation*}
 v^\epsilon = \mathsf{E}\left[ \int_0^T V_t\mathrm{d}t \right] = \int_0^T V_0(t)\mathrm{d}t,
\end{equation*}
as $\epsilon \to 0$, we have
\begin{equation*}
(X^\epsilon,Y^\epsilon) = 
\left(\frac{1}{\sqrt{v^\epsilon}}\int_t^T \sqrt{V_0(s)}\mathrm{d}B_s,
 \int_0^T V_0(s)\int_0^s\sum_{i=1}^d k_i(s-u) \mathrm{d}W^i_u \mathrm{d}s\right) + o_p(1).
\end{equation*}
The leading term $(X,Y)$ is centered normal
with covariance
\begin{equation*}
    \begin{split}
  \mathsf{E}[XY]   & =   \mathsf{E}\left[ 
        \frac{1}{\sqrt{v^\epsilon}}\int_0^T \sqrt{V_0(s)}\mathrm{d}B_s
 \int_0^T V_0(s)\int_0^s \sum_{i=1}^d k_i(s-u) \mathrm{d}W^i_u\mathrm{d}s\right]\\
        &= \frac{1}{\sqrt{v^\epsilon}} \sum_{i=1}^d
        \mathsf{E}\left[ 
        \int_0^T \sqrt{V_0(s)}\mathrm{d}B_s
 \int_0^T \int_s^T V_0(u)k_i(u-s)\mathrm{d}u  \mathrm{d}W^i_s
        \right]
        \\
        &= \frac{1}{\sqrt{v^\epsilon}}
        \int_0^T \sqrt{V_0(s)} \int_s^T V_0(u)\sum_{i=1}^d \rho_ik_i(u-s)\mathrm{d}u  \mathrm{d}s.
    \end{split}
\end{equation*}
The uniform integrability of $Y^\epsilon$ can be easily shown by observing its boundedness in $L^2$. Since $\mathsf{E}[Y|X=x] = \mathsf{E}[XY]x$, all the assumptions of Theorem~\ref{thm} are satisfied. We have in particular \eqref{skew} with $\mathsf{E}[XY]$ given above.

\section{Proofs}\label{proof}
\subsection{Proof of Theorem~\ref{thm}}
i). Take $\delta \in (0, 1-\rho^2)$ and let
\begin{equation*}
    \tau = \inf\left\{t \geq 0;  \rho^2 \int_0^t V^\epsilon_s \mathrm{d}s \geq (1-\delta) v^\epsilon\right\}.
\end{equation*}
Here we are going to show
\begin{equation}\label{fir}
    \mathsf{E}[f(S^\epsilon_T)]=\mathsf{E}[f(S^\epsilon_{T\wedge \tau})] + o(\epsilon).
\end{equation}
Since $f$ is bounded,
\begin{equation*}
       |\mathsf{E}[f(S^\epsilon_T)]-\mathsf{E}[f(S^\epsilon_{T\wedge \tau})]| \leq 2\|f\|_\infty \mathsf{P}[\tau \leq T].
\end{equation*}
Then, we obtain \eqref{fir} from
\begin{equation}\label{tau}
\begin{split}
    \mathsf{P}[\tau \leq T]&=\mathsf{P}[\epsilon \rho^2 Y^\epsilon \geq (1-\delta-\rho^2)v^\epsilon]\\ & \leq \frac{  \epsilon \rho^2}{(1-\delta-\rho^2)v^\epsilon}
    \mathsf{E}\left[|Y^\epsilon|; |Y^\epsilon| \geq 
    \frac{(1-\delta-\rho^2)v^\epsilon}{\epsilon \rho^2}
    \right],
\end{split}
\end{equation}
which is of $o(\epsilon)$
by the uniform integrability of $Y^\epsilon$ and \eqref{liminf}.
\\
ii). 
Here we decompose $ \mathsf{E}[f(S^\epsilon_{T\wedge \tau})]$ to extract its leading term.
Since $W^\perp$ is independent of $\mathscr{G}_T$, we have
\begin{equation*}
  \mathsf{E}[f(S^\epsilon_{T\wedge \tau})]   =\mathsf{E}[\mathsf{E}[f(S^\epsilon_{T\wedge \tau})| \mathscr{G}_T]] 
  = \mathsf{E}\left[
p\left(\hat{S}^\epsilon_T,(1-\rho^2) \int_0^{T\wedge \tau}V^\epsilon_s \mathrm{d}s \right)
        \right],
\end{equation*}
where
\begin{equation*}
    \hat{S}^\epsilon_t = S_0 \exp\left\{ 
\rho Z^\epsilon_t  - \frac{\rho^2}{2} \langle Z^\epsilon \rangle_t   \right\},\ \ Z^\epsilon_t = \int_0^{t \wedge \tau} \sqrt{V^\epsilon_s} \mathrm{d}W_s.
\end{equation*}
and
\begin{equation*}
    p(s,t) = \int f\left(s \exp\left\{\sqrt{t}x - \frac{t}{2} \right\}\right)\phi(x)\mathrm{d}x.
\end{equation*}
Since $p$ solves the partial differential equation
\begin{equation}\label{pde}
    \frac{\partial p}{\partial t} = \frac{1}{2} s^2 \frac{\partial^2 p}{\partial s^2}, \ \ 
    p(s,0) = f(s),
\end{equation}
putting
\begin{equation*}
\Sigma^\epsilon_u = v^\epsilon - \rho^2 \int_0^{u \wedge \tau} V^\epsilon_t \mathrm{d}t,
\end{equation*}
 It\^o's formula gives
\begin{equation*}
    p(\hat{S}^\epsilon_T,\Sigma^\epsilon_T) = p(S_0,v^\epsilon ) + \int_0^T\frac{\partial p}{\partial s}(\hat{S}^\epsilon_u,\Sigma^\epsilon_u) \mathrm{d}\hat{S}^\epsilon_u,
\end{equation*}
and so, noting $\Sigma^\epsilon_u \geq \delta v^\epsilon > 0$,
\begin{equation*}
    \mathsf{E}[ p(\hat{S}^\epsilon_T,\Sigma^\epsilon_T)] = p(S_0,v^\epsilon ).
\end{equation*}
Therefore,
\begin{equation}\label{sec}
\begin{split}
&    \mathsf{E}[f(S^\epsilon_{T\wedge \tau})]\\ 
&= 
 p(S_0,v^\epsilon ) + \mathsf{E}\left[
p\left(\hat{S}^\epsilon_T, (1-\rho^2) \int_0^{T\wedge \tau}V^\epsilon_u \mathrm{d}u\right)
-p(\hat{S}^\epsilon_T,\Sigma^\epsilon_T)
        \right]\\
&= 
 p(S_0,v^\epsilon ) + \mathsf{E}\left[
p\left(\hat{S}^\epsilon_T,\Sigma^\epsilon_T + \epsilon \hat{Y}^\epsilon\right)
-p(\hat{S}^\epsilon_T,\Sigma^\epsilon_T)
        \right],
        \end{split}
\end{equation}
where
\begin{equation*}
    \hat{Y}^\epsilon = \frac{1}{\epsilon}\left(
\int_0^{T\wedge \tau}V^\epsilon_u \mathrm{d}u -v^\epsilon
    \right) = Y^\epsilon 1_{\{\tau > T\}} + \frac{(1-\delta-\rho^2)v^\epsilon}{\epsilon \rho^2} 1_{\{\tau \leq T\}}.
\end{equation*}
iii). Here we are going to show that $\hat{Y}^\epsilon$ is uniformly integrable.
By \eqref{tau},
\begin{equation*}
\begin{split}
& \sup_{\epsilon > 0}    \mathsf{E}\left[
     \frac{(1-\delta-\rho^2)v^\epsilon}{\epsilon \rho^2} 1_{\{\tau \leq T\}} ;
      \frac{(1-\delta-\rho^2)v^\epsilon}{\epsilon \rho^2} 1_{\{\tau \leq T\}} \geq K\right] \\
      & \leq \sup\left\{
       \frac{(1-\delta-\rho^2)v^\epsilon}{\epsilon \rho^2} \mathsf{P}[\tau \leq T]; 
      \frac{(1-\delta-\rho^2)v^\epsilon}{\epsilon \rho^2} \geq K
     \right\}
      \\
      & \leq\sup\left\{
       \mathsf{E}\left[|Y^\epsilon|; |Y^\epsilon| \geq 
    \frac{(1-\delta-\rho^2)v^\epsilon}{\epsilon \rho^2}
    \right];  
     \frac{(1-\delta-\rho^2)v^\epsilon}{\epsilon \rho^2} \geq K
    \right\}
    \\
    & \leq \sup_{\epsilon > 0}\mathsf{E}\left[|Y^\epsilon|; |Y^\epsilon| \geq  K \right],
      \end{split}
\end{equation*}
which converges to $0$ as $K\to \infty$ by the uniform integrability of $Y^\epsilon$. Therefore, $\hat{Y}^\epsilon$ also is uniformly integrable.\\
iv). Let
\begin{equation*}
    A = \left\{ (1-\rho^2)\int_0^{T\wedge \tau} V^\epsilon_u \mathrm{d}u > \delta v^\epsilon \right\}.
\end{equation*}
We have
\begin{equation*}
\begin{split}
    \mathsf{P}[A^c] &= \mathsf{P}\left[
    \hat{Y}^\epsilon \leq - \frac{(1-\delta - \rho^2)v^\epsilon}{\epsilon (1-\rho^2)}
    \right]
    \\ & \leq \frac{\epsilon (1-\rho^2)}{(1-\delta - \rho^2)v^\epsilon} \mathsf{E}\left[|\hat{Y}^\epsilon|; |\hat{Y}^\epsilon| \geq \frac{(1-\delta - \rho^2)v^\epsilon}{\epsilon (1-\rho^2)} \right],
\end{split}
\end{equation*}
which is of $o(\epsilon)$ by the uniform integrability of $\hat{Y}^\epsilon$ and \eqref{liminf}.
Since $f$ is bounded, so is $p$. This implies that
\begin{equation*}
     \mathsf{E}\left[
p\left(\hat{S}^\epsilon_T,\Sigma^\epsilon_T + \epsilon \hat{Y}^\epsilon\right)
-p(\hat{S}^\epsilon_T,\Sigma^\epsilon_T)
        \right]
        =  \mathsf{E}\left[
p\left(\hat{S}^\epsilon_T,\Sigma^\epsilon_T + \epsilon \hat{Y}^\epsilon\right)
-p(\hat{S}^\epsilon_T,\Sigma^\epsilon_T) ; A
        \right] + o(\epsilon)
\end{equation*}
We have
\begin{equation*}
    p\left(\hat{S}^\epsilon_T,\Sigma^\epsilon_T + \epsilon \hat{Y}^\epsilon\right)
-p(\hat{S}^\epsilon_T,\Sigma^\epsilon_T)
= \epsilon \hat{Y}^\epsilon \int_0^1 \frac{\partial p}{\partial t}
\left(\hat{S}^\epsilon_T,\Sigma^\epsilon_T + r \epsilon \hat{Y}^\epsilon\right) \mathrm{d}r 
\end{equation*}
and on the set $A$,
\begin{equation*}
    \Sigma^\epsilon_T + r \epsilon \hat{Y}^\epsilon \geq \delta v^\epsilon
\end{equation*}
for all $ r \in [0,1]$. Combining with the above, we obtain
\begin{equation}\label{thi}
\begin{split}
&     \mathsf{E}\left[
p\left(\hat{S}^\epsilon_T,\Sigma^\epsilon_T + \epsilon \hat{Y}^\epsilon\right)
-p(\hat{S}^\epsilon_T,\Sigma^\epsilon_T)
        \right]
       \\ & =   \epsilon  \mathsf{E}\left[\hat{Y}^\epsilon
\int_0^1 \frac{\partial p}{\partial t}
\left(\hat{S}^\epsilon_T,(\Sigma^\epsilon_T + r \epsilon \hat{Y}^\epsilon)\vee (\delta v^\epsilon)\right) \mathrm{d}r; A
        \right] + o(\epsilon)
\end{split}
\end{equation}
v). Notice that
\begin{equation}
    \frac{\partial p}{\partial t}(s,t) = \int f(se^w)\frac{\partial}{\partial t} \phi\left(
\frac{w}{\sqrt{t}} + \frac{\sqrt{t}}{2}\right) \mathrm{d}w
\end{equation}
and because $f$ is bounded,
\begin{equation*}
     \frac{\partial p}{\partial t}(s,u\vee(\delta v^\epsilon))
\end{equation*}
is bounded in $s$, $u$ and sufficiently small $\epsilon$ under \eqref{liminf}.\\
vi). 
Since $(X^\epsilon, Y^\epsilon) \to (X,Y)$ in law, with the aid of i), we have
\begin{equation*}
\left( \frac{1}{\sqrt{v^\epsilon}}\int_0^{T\wedge \tau}\sqrt{V^\epsilon_t}\mathrm{d}B_t, \hat{Y}^\epsilon \right) \to 
\left( X,Y\right)
\end{equation*}
in law as $\epsilon \to 0$. On the other hand,
\begin{equation*}
    \int_0^{T\wedge \tau}\sqrt{V^\epsilon_t}\mathrm{d}B_t = \rho Z^\epsilon_T + \sqrt{(1-\rho^2)(v^\epsilon  + \epsilon \hat{Y}^\epsilon}) N^\epsilon,
\end{equation*}
where $N^\epsilon$ is a standard normal random variable independent of $(Z^\epsilon_T,\hat{Y}^\epsilon)$.
Therefore, the sequence of joint distribution
\begin{equation*}
    \left( \frac{1}{\sqrt{v^\epsilon}}\int_0^{T\wedge \tau}\sqrt{V^\epsilon_t}\mathrm{d}B_t, \hat{Y}^\epsilon, Z^\epsilon_T, N^\epsilon \right)
\end{equation*}
is tight, of which an accumulation point is uniquely determined as $(X,Y,Z,N)$
such that $N$ is a standard normal random variable independent of $(Y,Z)$ and
\begin{equation*}
    \rho Z + \sqrt{(1-\rho^2)v^0}N = X.
\end{equation*}
Since $X$ must be a standard normal random variable by the martingale central limit theorem,
$Z$ is a centered normal random variable with variance $v^0$.\\
vii). From iii), iv), v) and vi), we have
\begin{equation}\label{fou}
\begin{split}
       & \lim_{\epsilon \to 0} \mathsf{E}\left[\hat{Y}^\epsilon
\int_0^1 \frac{\partial p}{\partial t}
\left(\hat{S}^\epsilon_T,(\Sigma^\epsilon_T + r \epsilon \hat{Y}^\epsilon)\vee (\delta v^\epsilon)\right) \mathrm{d}r; A
        \right]  \\ & = \lim_{\epsilon \to 0} \mathsf{E}\left[\hat{Y}^\epsilon
\int_0^1 \frac{\partial p}{\partial t}
\left(\hat{S}^\epsilon_T,(\Sigma^\epsilon_T + r \epsilon \hat{Y}^\epsilon)\vee (\delta v^\epsilon)\right) \mathrm{d}r
        \right]\\
        &= \mathsf{E}\left[ Y \frac{\partial p}{\partial t}(\hat{S}^0,(1-\rho^2)v^0)\right],
\end{split}
\end{equation}
where
\begin{equation*}
    \hat{S}^0 = S_0 \exp\left\{\rho Z - \frac{\rho^2}{2}v^0\right\}.
\end{equation*}
From \eqref{fir}, \eqref{sec}, \eqref{thi} and \eqref{fou}, we have
\begin{equation}\label{fiv}
    \mathsf{E}[f(S^\epsilon_T)] = p(S_0,v^\epsilon) + \epsilon \mathsf{E}\left[ Y \frac{\partial p}{\partial t}(\hat{S}^0,(1-\rho^2)v^0)\right] + o(\epsilon).
\end{equation}
viii). Let $\check{W}$ be a standard Brownian motion independent of $(Y,Z)$ and 
\begin{equation*}
    \check{S}_\tau = \hat{S}^0\exp\left\{\sqrt{(1-\rho^2)v^0}\check{W}_\tau - \frac{1}{2}(1-\rho^2)v^0 \tau \right\}, \ \ \check{\Sigma}_\tau = (1-\rho^2)v^0(1-\tau) + \delta.
\end{equation*}
Then, by It\^o's formula and \eqref{pde},
\begin{equation*}
    \frac{\partial p}{\partial t}(\check{S}_1,\check{\Sigma}_1) = 
    \frac{\partial p}{\partial t}(\check{S}_0,\check{\Sigma}_0) + \int_0^1 
    \frac{\partial^2 p}{\partial s \partial t}(\check{S}_\tau,\check{\Sigma}_\tau)\mathrm{d}\check{S}_\tau,
\end{equation*}
which implies
\begin{equation}\label{six}
    \mathsf{E}\left[ Y \frac{\partial p}{\partial t}(\check{S}_1,\check{\Sigma}_1)
    \right] = 
    \mathsf{E}\left[ Y \frac{\partial p}{\partial t}(\check{S}_0,\check{\Sigma}_0)
    \right] = 
    \mathsf{E}\left[ Y \frac{\partial p}{\partial t}(\hat{S}^0,(1-\rho^2)v^0 + \delta)\right].
\end{equation}
Notice that the joint law of $(\check{S}_1,Y)$ is identical to
\begin{equation*}
   \left( S_0 \exp\left\{ \sqrt{v^0}X - \frac{v^0}{2} \right\}, Y \right).
\end{equation*}
Therefore, using again \eqref{pde},
\begin{equation*}
    \begin{split}
     &   \mathsf{E}\left[ Y \frac{\partial p}{\partial t}(\check{S}_1,\check{\Sigma}_1)
    \right]
    = \frac{1}{2} \mathsf{E}\left[ Y \hat{S}_1^2\frac{\partial^2 p}{\partial s^2}(\check{S}_1,\check{\Sigma}_1)
    \right]\\
    &= \frac{1}{2}\int
    \mathsf{E}[Y|X=x] \left(S_0e^{\sqrt{v^0}x -v^0/2}\right)^2
    \frac{\partial^2 p}{\partial s^2}\left(S_0e^{\sqrt{v^0}x -v^0/2},\delta\right) \phi(x)
    \mathrm{d}x\\
    &= \frac{1}{2}\int
    \mathsf{E}[Y|X=x] 
    \left(\frac{1}{v^0}
    \frac{\partial^2 }{\partial x^2}
    -  \frac{1}{\sqrt{v^0}}
    \frac{\partial }{\partial x}\right)p\left(S_0e^{\sqrt{v^0}x -v^0/2},\delta\right)\phi(x)\mathrm{d}x
    \\
    &=
    \frac{1}{2}\int
    p\left(S_0e^{\sqrt{v^0}x -v^0/2},\delta\right)
    \left(\frac{1}{v^0}
    \frac{\partial^2 }{\partial x^2}
    +  \frac{1}{\sqrt{v^0}}
    \frac{\partial }{\partial x}\right) (\mathsf{E}[Y|X=x] \phi(x))\mathrm{d}x.
    \end{split}
\end{equation*}
Here we have used \eqref{boundary}. Then, by \eqref{six},
\begin{equation*}
\begin{split}
     &\mathsf{E}\left[ Y \frac{\partial p}{\partial t}(\hat{S}^0,(1-\rho^2)v^0)\right]
     = \lim_{\delta \to 0}
      \mathsf{E}\left[ Y \frac{\partial p}{\partial t}(\hat{S}^0,(1-\rho^2)v^0 + \delta)\right]
      \\& =
      \lim_{\delta \to 0}
       \frac{1}{2}\int
    p\left(S_0e^{\sqrt{v^0}x -v^0/2},\delta\right)
    \left(\frac{1}{v^0}
    \frac{\partial^2 }{\partial x^2}
    +  \frac{1}{\sqrt{v^0}}
    \frac{\partial }{\partial x}\right) (\mathsf{E}[Y|X=x] \phi(x))\mathrm{d}x
      \\& = \frac{1}{2}\int
    p\left(S_0e^{\sqrt{v^0}x -v^0/2},0\right)
    \left(\frac{1}{v^0}
    \frac{\partial^2 }{\partial x^2}
    +  \frac{1}{\sqrt{v^0}}
    \frac{\partial }{\partial x}\right) (\mathsf{E}[Y|X=x] \phi(x))\mathrm{d}x.
\end{split}
\end{equation*}
Here, the last equality follows from the Fubini-Tonelli theorem and that
\begin{equation*}
    z \mapsto 
    \int
    f\left(S_0e^{\sqrt{v^0}x -v^0/2 + z}\right)
    \left(\frac{1}{v^0}
    \frac{\partial^2 }{\partial x^2}
    +  \frac{1}{\sqrt{v^0}}
    \frac{\partial }{\partial x}\right) (\mathsf{E}[Y|X=x] \phi(x))\mathrm{d}x
\end{equation*}
is bounded and continuous, which can be shown by
the generalized dominated convergence theorem.
The rest would be clear.
 \hfill{$\square$}
 \subsection{Proof of Corollary~\ref{cor1}}
First note that
\begin{equation}
    p_K(S_0,v^\epsilon ) =  \int f\left(S_0 \exp\left\{\sqrt{v^\epsilon}x - \frac{v^\epsilon}{2} \right\}\right)\phi(x)\mathrm{d}x
\end{equation}
for $f(s) = (K-s)_+$. Let $g(x) = \mathsf{E}[Y|X = x] \phi(x)$ 
and $d_- = d_-(S_0,v^\epsilon )$. We have
\begin{equation}
    \begin{split}
         &\int f\left(S_0 \exp\left\{\sqrt{v^\epsilon}x - \frac{v^\epsilon}{2} \right\}\right)\frac{\mathrm{d}}{\mathrm{d}x} g(x)\mathrm{d}x 
         \\ & = \int_{-\infty}^{-d_-}
       \left(K-S_0 \exp\left\{\sqrt{v^\epsilon}x - \frac{v^\epsilon}{2} \right\}\right)
       \frac{\mathrm{d}}{\mathrm{d}x} g(x)\mathrm{d}x \\
       &= \sqrt{v^\epsilon}\int_{-\infty}^{-d_-}
       S_0 \exp\left\{\sqrt{v^\epsilon}x - \frac{v^\epsilon}{2} \right\}
       g(x)\mathrm{d}x
    \end{split}
\end{equation}
by \eqref{boundary}
and similarly,
\begin{equation}
    \begin{split}
         &\int f\left(S_0 \exp\left\{\sqrt{v^\epsilon}x - \frac{v^\epsilon}{2} \right\}\right)\frac{\mathrm{d}^2}{\mathrm{d}x^2} g(x)\mathrm{d}x \\
       &= \sqrt{v^\epsilon}\int_{-\infty}^{-d_-}
       S_0 \exp\left\{\sqrt{v^\epsilon}x - \frac{v^\epsilon}{2} \right\}
       \frac{\mathrm{d}}{\mathrm{d}x} g(x)\mathrm{d}x \\
       &= \sqrt{v^\epsilon}\left[ 
       S_0 \exp\left\{\sqrt{v^\epsilon}x - \frac{v^\epsilon}{2} \right\} g(x)
       \right]^{-d_-}_{-\infty} - v^\epsilon\int_{-\infty}^{-d_-}
       S_0 \exp\left\{\sqrt{v^\epsilon}x - \frac{v^\epsilon}{2} \right\}
       g(x)\mathrm{d}x.
    \end{split}
\end{equation}
The result then follows from Theorem~\ref{thm}.
\hfill{$\square$}

\end{document}